\documentclass[11pt]{amsart} 

\renewcommand{\bar}{\overline}
\newcommand{\eps}{\epsilon}
\newcommand{\pa}{\partial}
\newcommand{\ph}{\varphi}

\newcommand{\rr}{{r \leq\frac{\log m}{\sqrt m}}}  

\renewcommand{\eps}{\varepsilon} 

\font\strange=msbm10

\newcommand{\C}{{{\mathchoice  {\hbox{$\textstyle{\text{\strange C}}$}}
{\hbox{$\textstyle{\text{\strange C}}$}}
{\hbox{$\scriptstyle{\text{\strange C}}$}}
{\hbox{$\scriptscriptstyle{\text{\strange C}}$}}}}}

\newcommand{\R}{{{\mathchoice  {\hbox{$\textstyle{\text{\strange R}}$}}
{\hbox{$\textstyle{\text{\strange R}}$}}
{\hbox{$\scriptstyle  N\kern-0.3em  R$}}  
{\hbox{$\scriptscriptstyle  R\kern-0.2em  R$}}}}}

\newcommand{\Z}{{{\mathchoice  {\hbox{$\textstyle{\text{\strange Z}}$}}
{\hbox{$\textstyle{\text{\strange Z}}$}}
{\hbox{$\scriptstyle  Z\kern-0.3em  Z$}}
{\hbox{$\scriptscriptstyle  Z\kern-0.2em  Z$}}}}}

\newcommand{\N}{{{\mathchoice  {\hbox{$\textstyle{\text{\strange N}}$}}
{\hbox{$\textstyle{\text{\strange N}}$}}
{\hbox{$\scriptstyle  N\kern-0.3em  N$}}
{\hbox{$\scriptscriptstyle  N\kern-0.2em  N$}}}}}

\newcommand{\frk}[1]{{\mathfrak{#1}}}

\newcommand{\oo}[1]{{O(\frac{1}{m^{{#1}}})}}
\newcommand{\ooo}[1]{{O(|z|^{{#1}})}}
\newcommand{\bb}{{\frac{\sqrt{-1}}{2\pi}}}

\newcommand{\crr}[4]{{R_{{#1}\bar{#2}{#3}\bar{#4}}}}
\newcommand{\css}[2]{{R_{{#1}\bar{#2}}}}

\newcommand{\cp}{{\mathbb C}P^n}
\newcommand{\ka}{K\"ahler }
\renewcommand{\phi}{\varphi}
\usepackage{amsmath,amsthm,amscd}

\usepackage{calc}               
  
\title
[]{The log term of  Szeg\"o Kernel}
\date{October 24, 2001}


\author{Zhiqin Lu} 
\address[Zhiqin Lu] 
{Department of Mathematics\\
University of California at Irvine\\
Irvine, CA 92612}
\email{zlu@math.uci.edu}
\thanks{The first author is supported by  NSF grant DMS
0204667 and the Alfred P. Sloan Fellowship. The second author
is supported by an NSF grant.}

\author{Gang Tian}
\address[Gang Tian]
{Department of Mathematics\\
Massachusetts Institute of Technology\\
Cambridge, MA 02139}
\email{tian@math.mit.edu} 

\newtheorem{theorem}{Theorem}[section]

\newtheorem*{qu}{Conjecture}
\newtheorem{lemma}{Lemma}[section]
\newtheorem{cor}{Corollary}[section]
\newtheorem{prop}{Proposition}[section] 

\newtheorem{definition}{Definition}[section]

\theoremstyle{remark}
\newtheorem{rem}{Remark}[section] 

\begin{document}
\maketitle

\numberwithin{equation}{section}

\tableofcontents 

\section{Introductions}
%2003.2.18.9.22

Prescribing geometric structures of a complex manifold
often introduces interesting and important 
partial differential equations. 
A typical example of this kind is the problem 
of finding the \ka metrics 
with constant scalar curvature
on a \ka manifold. 
Such a problem defines 
a  fourth order elliptic 
partial differential equation. The study of these
partial differential equations, including the
K\"ahler-Einstein equations, forms one of the richest
topics in complex geometry.

In this paper, we  introduce a new set of 
equations coming from the Szeg\"o kernel (Bergman
kernel, resp.) of a  unit circle (unit disk, resp)
bundle. We prove that these
equations, which generalize the equation of 
finding
\ka
 metrics with constant scalar curvature,
are all elliptic. As an application of the result, we
relate the Ramadanov Conjecture to these equations
and prove a local rigidity theorem concerning the log
term of the Szeg\"o kernel.

Our basic setting is as follows:
let $(L,h)\rightarrow M$ be a positive Hermitian line
bundle over the compact complex  manifold $M$ of
dimension $n$. 
The pair $(M,L)$ is called a polarised manifold. 
The \ka metric $\omega$ of $M$ is defined to be the 
curvature of the Hermitian metric $h$.
Let $L^*$ be the dual bundle of $L$.  The
unit circle bundle $X$
 of $L^*$ is a strictly pseudoconvex
manifold, with the natural measure defined by the $S^1$
action and the polarization of $M$. That is, the measure
is
$dV=\frac{1}{n!}\pi^*(\omega^n)\wedge d\theta$, where
$\frac{\pa}{\pa\theta}$ is the infinitesimal $S^1$ action
on the unit circle bundle.
 The
Szeg\"o projection
$\Pi$ is a linear map from
$L^2(X)$ to the Hardy space $H^2(X)$, which is the 
space of $L^2$ boundary
functions of holomorphic functions of the unit disk
bundle $D$.
 Let $\Pi(x,y)$ be the
Szeg\"o kernel of $X$, i.e.,   
$\Pi(x,y)$ is the function on $X\times X$ such that 
for any $f\in L^2(X)$, 
$\int_X\Pi(x,y)f(y)dy\in H^2(X)$, where $dy=dV$ is the
measure
defined above. Then by ~\cite{BS1}, there is
a paramatrix
\[
s(x,y,t)\sim\sum_{k=0}^\infty t^{n-k} s_k(x,y),
\]
where $s_k(x,y)\,(k\in\Z_+)$ are smooth functions on $X\times X$ and
$t\in \R$,
such that
\begin{equation}\label{1-9}
\Pi(x,y)=\int_0^\infty e^{it\psi(x,y)} s(x,y,t) dt
\end{equation}
for some suitable complex phase function $\psi(x,y)$
of $X\times X$. 

In general, the paramatrix of the Szeg\"o kernel 
of a pseudoconvex manifold
is
quite difficult to compute. However, since the bundle
$X$
is $S^1$ invariant, we 
may split the Szeg\"o kernel into  several pieces.
More precisely, Let
$\frac{\pa}{\pa\theta}$ be the infinitesimal $S^1$ action 
of $X$. Define
\[
H^2_m(X)=\{f\in H^2(X)| \frac{\pa}{\pa \theta}
f=\sqrt{-1}m f\}.
\]
Let $\Pi_m$ be the projection of $H^2(X)$ to $H^2_m(X)$.
Then the kernel $\Pi_m(x,y)$ of $L^2(X)\rightarrow
H^2_m(X)$ 
is the Fourier coefficient of $\Pi(x,y)$:
\[
\Pi_m(x,y)=\frac{1}{2\pi}\int_{S^1}
\Pi(x,r_\theta
y)e^{m\sqrt{-1}\theta}d\theta.
\]
Using the paramatrix of the Szeg\"o kernel,
Zelditch (and Catlin~\cite{Cat} independently for the
Bergman kernel) was able to prove that there is an
asymptotic expansion of
$\Pi_m(x,x)$ (cf. Theorem~\ref{zelditch})
\begin{equation}\label{kkk}
\Pi_m(x,x)\sim m^n(a_0+\frac{a_1}{m}+\cdots),
\end{equation}
where $a_k$'s are all smooth functions of $M$. 
The expansion is called Tian-Yau-Zelditch expansion.
In~\cite{Lu10}, the first author was able to prove that
all $a_k$'s are polynomials of the curvature
and its derivatives. In particular, $a_0=1$
and $a_1=\rho$, the scalar curvature of the \ka manifold.
Thus the equation of finding the metrics such that
$a_1={\rm const}$ is the equation of finding the
\ka metrics with
constant scalar curvature.

Because of the work of Donaldson~\cite{SKD1}, 
it is  natural to study metrics with
$a_k$ being prescribed for $k\geq 2$.
Donaldson was interested in modifying $h^m$ to a Hermitian
metric $h'$ for some large $m$ such that the metric $h'$
is balanced. As a corollary of his result, Donaldson
was able to give  a proof of the uniqueness of 
the \ka metrics
of constant scalar curvature. Since $a_1=\frac 12\rho$, where 
$\rho$ is the scalar curvature, $\int_M(a_1-\bar{a_1})
\theta$ defines the Futaki invariants, where $\bar{a_1}$
is the average of $a_1$ and $\theta$ is the Hamiltonian
function of a holomorphic vector field. Nonlinearizing
the Futaki invariants we get the Mabuchi's $K$ energy,
whose convexity plays the key role in proving the uniqueness
of the metrics of constant sclar curvature (cf.
~\cite{chen1}).

We wish to study the analogue problems for $a_k$ when
$k>1$.
 In this
paper, among the other results,
we prove that for any given $k$ and function $f$,
for a fixed metric $\omega$, the equation of finding
the function $\ph$ such that 
$a_k(\omega-\bb\pa\bar\pa\ph)=f$ is an 
elliptic equation of order $2k+2$. Thus prescribing $a_k$
gives an interesting set of new elliptic equations.

Since the Bergman potential 
$\Pi_m(x,x)$ being a constant
implies stability(cf. ~\cite{Luo}, ~\cite{Zh}),
we are particularly interested in the question of
finding metrics such that $a_k=0$ for $k>n$.
Such a question is related to the Ramadanov
Conjecture~\cite{rama}. The conjecture,
in terms of the Bergman kernel, can be stated as
follows:  

\begin{qu}[Ramadanov~\cite{rama}]
Let $\Omega$ be a bounded strongly pseudoconvex domain
of $\C^n$. Assume that the log term of the
Bergman kernel is zero, then $\Omega$ is
biholomorphic equivalent to the unit ball of 
$\C^n$.  
\end{qu}  
 
Not much is
known about the Conjecture for $n>1$. If $\Omega$ is a 
complete Rienhardt domain
of $\C^2$, the conjecture was proved by
Nakazawa~\cite{Naka}.
  The
conjecture was proved to be true for any
strongly pseudoconvex domain in $\C^2$ by 
Graham~\cite{G2}, using an unpublished note
of Burns. In~\cite{HKN}, the computation needed in Graham's
proof was given. Boutet de Monvel~\cite{Boutet} gave an
independent proof of Graham's result around the same time.

There are only partial results in higher
dimensions. K. Hirachi~\cite{HK3}
proved that the Radamanov conjecture is true for real
ellipsoids  that are
sufficiently close to the unit ball. In the Szeg\"o
kernel case, he~\cite{HK4} proved that if $n=2$ with
tranversal symmetry and if the  log term of the Szeg\"o kernel
vanishes to the third order, then the boundary is spherical.
Hanges~\cite{Hanges} proved a similar result with additional
assumption on the choice of volume element on the boundary.
See~\cite{HK3} for further references of the Conjecture.

One can form the similar conjecture for the Szeg\"o
kernel as well. The same conjecture makes
sense if we replace the bounded
strongly pseudoconvex domains by strongly 
pseudoconvex manifolds. 

Let $H^*$ be the universal line bundle of the 
complex projective space $\cp$. The unit
circle bundle $X$ of $H^*$ is the unit sphere in $\C^n$. 
With this observation, we form the following 
Ramadanov Conjecture for the unit circle bundle $X$:

\begin{qu}
Let $\omega\in[\omega_{FS}]$ be a \ka metric on $\cp$ 
which is in the same cohomology class as the 
Fubini-Study metric $\omega_{FS}$.
Let $(H,h)$ be the hyperplane bundle whose 
curvature
is $\omega$. Let $X$ be the unit circle bundle
of the universal line bundle $H^*$. 
If the log term of the Szeg\"o kernel of $X$ vanishes, 
then there
is an automorphism
$\ph:\cp\rightarrow \cp$ such that
$\ph^*\omega=\omega_{FS}$.
\end{qu}

We confirmed the above conjecture for $n=1$
in this paper.
In the case $n=1$, the unit circle bundle 
$X$ is of dimension
$3$ and the result is parallel to the case of 
strongly pseudoconvex domains
in $\C^2$ in the Ramadanov Conjecture. 

\begin{theorem}[The case $n=1$]\label{fud2}
Let $\omega\in[\omega_{FS}]$ be a \ka form on $\C P^1$ 
which  is in the same cohomology class as the 
Fubini Study metric $\omega_{FS}$.
Let $(H,h)$ be the hyperplane bundle whose 
curvature
is $\omega$. Let $X$ be the unit circle bundle
of the universal bundle $H^*$. 
If the log term of the Szeg\"o kernel of $X$ vanishes, 
then there
is an automorphism
$\ph:\C P^1\rightarrow \C P^1$ such that
$\ph^*\omega=\omega_{FS}$.
\end{theorem}

The Conjecture is still open in high dimensions.
The main result of this paper is the following
local rigidity theorem:

\begin{theorem}\label{fud1}
Let $h$ be the standard metric on $(\cp, H)$. Let $h'$
be another metric on the hyperplane bundle $H$ over
$\cp$.
Then there is an $\eps>0$, depending only on $n$, such
that for any $h'$ with
\[
||h'/h-1||_{C^{2n+4}}<\eps,
\]
and  the log term of the Szeg\"o kernel of the unit
circle bundle of
$h'$ being zero, there is an automorphism $f$  of $\cp$
such that
$f^*(\omega_{h'}) =\omega_{FS}$, where $\omega_{h'}$ is
the curvature of
$h'$.
\end{theorem}

The organization of this paper is as follows:
In \S 2,
we prove that if the log term is zero, then 
$a_k=0$ for $k>n$. This result let us
use the methods of partial differential equations to
study the Ramadanov Conjecture. In
\S 3, by tracing the terms in $a_k$ of the highest 
Weyl weight,  we prove
that the equations 
$a_k=f$
are all elliptic
equations. As a corollary, we get the Schauder estimate
(Corollary~\ref{cor31}).

In Theorem~\ref{thm52} and Theorem~\ref{34}, we study
the uniformity of the Tian-Yau-Zelditch expansion.
Using the elliptic estimate, we have the following
general result:

\begin{theorem}\label{fud}
If $h$ is the metric such that the log term of the 
Szeg\"o kernel (of the unit circle bundle $X$ 
of $L^*$, the dual bundle of the ample line bundle
$L$ over $M$)
vanishes, then there is a finite 
dimensional
vector subspace $V$ 
of the space of smooth functions on $M$
such that if $\ph\notin V$, the log term of the Szeg\"o kernel
is not zero
for the metric $he^{\eps\ph}$ for sufficiently small 
$\eps$.
\end{theorem}

This proves that generically, the log term is not zero.

The technical heart of this paper is in \S6 and \S 7.
In
\S 6,  we computed
concretely the vector space $V$ in Theorem~\ref{fud}.
In order to determine the vector space $V$ in
Theorem~\ref{fud}, we use the fact that the orthonormal 
basis of $H^0(\cp, H^m)$ can always be  explicitly
written. Then Lemma~\ref{lem61} and
Proposition~\ref{comb} become purely combinatoric.
From these results, we prove that $V$ must be the
eigenspace of some eigenvalue of $\cp$.
In Theorem~\ref{lamb}, we further confirmed that the
eigenvalue must be $(n+1)$. It is well known that 
the eigenfunctions with respect to the eigenvalue
$(n+1)$ of $\cp$
are the Hamiltonian  functions of holomorphic
vector fields on $\cp$. Thus in order to prove Theorem
~\ref{fud1}, we have to get rid of the actions of the
automorphism group of $\cp$ generated by the holomorphic
vector fields (the idea was first used by Bando-Mabuchi
~\cite{BM}).
In ~\cite{CT2}, the authors introduced the concept of a
\ka metric being centrally positioned. In our case, we
need to know that if a metric is not far away from the
Fubini-Study metric, then whether we can find a ``small''
automorphism under which the metric is centrally
positioned. This is done in Lemma~\ref{lem72} using the
contraction principal. 

{\bf Acknowledgement.} We thank K. Hirachi
for the references of the history of the 
Ramadanov Conjecture.

\section{Pseudoconvex manifolds with zero log term}
%2003.02.18.9.23

In this section, we study the relations
between the vanishing of the log term in the  Szeg\"o
or Bergman kernels and the coefficients in the 
Tian-Yau-Zelditch expansion~\cite{Lu10}. Except for the
last theorem, or otherwise stated,  most results of
this section were known by the work of
Zelditch~\cite{sz} and Catlin~\cite{Cat}.
We begin by the following standard settings.

Suppose $(L,h)\rightarrow M$ is a positive Hermitian line bundle over the
compact complex manifold $M$
as in the previous section.
Let $(L^*,h^{-1})$ be the dual bundle. 
 We define a
smooth function $\rho :L^*\rightarrow \R$ as follows:
let $U\subset M$ be an open neighborhood of $M$ such
that $L^*|_U\overset{\ph}{=}U\times \C$ is a local
trivialization. Let
\[
\rho(x)=\frac{1}{h(z)}|v|^2-1,
\]
where $h(z)$ is the local representation of the
Hermitian metric $h$ under the trivialization $\ph$ and
$\ph(x)=(z,v)$. It is not hard to see that 
$\rho(x)$ does not depend on the choice of the local trivialization.

\begin{definition}
Let $D=\{x\in L^*|\rho(x)\leq 0\}$. Let $X$ be the boundary $\pa D$
of $D$. We
call $D$ and $X$ the unit disk bundle and the unit circle bundle of $L^*$,
respectively. 
\end{definition}

In what follows, we will take the Szeg\"o
kernel as an example in our proof. The
results of the Bergman kernel are similar
and we will only  state the theorem without proof.

By definition, the curvature of $h$ is positive
and thus
$X$ is a strongly pseudoconvex manifold.

$X$ is $S^1$ invariant. 
Let $r_\theta: X\rightarrow X$ 
be defined by $r_\theta (z,v)=(z, v
e^{\sqrt{-1}\theta})$.
Let $\frac{\pa}{\pa \theta}$ be the infinitesimal
action of $S^1$ on $X$ defined by $r_\theta$. Let $\pi:
L^*\rightarrow M$ be the projection. Define the
measure 
$d\mu=\frac{1}{n!}\pi^*\omega_g^n\wedge d\theta$ on $X$.
The Szeg\"o kernel is defined as the kernel
 of the projection of the
space   $L^2(X)$ to the space $H^2(X)$, 
the Hardy space. By
definition, $H^2(X)$  is the space of $L^2$ functions
on 
$X$ which are the boundary
values of holomorphic functions on $D$. 

Let $\Pi$ be the Szeg\"o projection and 
let $\Pi(x,y)$
be its kernel. Then
\[
\Pi: L^2(X)\rightarrow H^2(X)
\]
such that for any $f\in L^2(X)$,
\[
\Pi f=\int_X\Pi(x,y)f(y) d\mu(y)
\]
is  in $L^2(X)$ and is the boundary function of a holomorphic function on
$D$.

Suppose $M$ is covered by finite coordinate
charts $\{U_\alpha\}_{\alpha\in I}$. The transition
functions of the line bundle $L$ are a set of
holomorphic functions $g_{\alpha\beta}$ on
$U_\alpha\cap U_\beta\neq\emptyset$. Let
$h_\alpha(x)$ be the local representation
of the Hermitian metric $h$. Then on 
$U_\alpha\cap U_\beta\neq\emptyset$ we have
\[
h_\alpha|g_{\alpha\beta}|^2=h_\beta.
\]
For each $h_\alpha$, we define a $C^\infty$ function 
$\tilde h(x,y)$ on $U_\alpha\times U_\alpha$ such that
$\tilde h_\alpha(x,x)=h_\alpha(x)$ and
$\tilde h_\alpha(x,y)$ is almost analytic in the sense
that
$\bar\pa_x \tilde h_\alpha(x,y)$ and $\pa_y\tilde
h(x,y)$ vanish at $x=y$ to infinite order. Such a
function
$\tilde h_\alpha(x,y)$ exists by~\cite{BS1}.

Let $\sigma_\alpha$ is the partition of the unity
subordinated to the covering $\cup \{U_\alpha\}$
such that $\sqrt{\sigma_\alpha}$ are all smooth.
 Define
\[
h_\alpha(x,y)=\sum_{U_\gamma\cap
U_\alpha\neq\emptyset}\sqrt{\sigma_\gamma(x)}
\sqrt{\sigma_\gamma(y)}
g_{\gamma\alpha}(x)\bar{g_{\gamma\alpha}(y)}
\tilde h_\gamma(x,y).
\]
It is then easy to check that
\[
h_\alpha(x,y)=h_\beta(x,y) g_{\beta\alpha}(x)
\bar{g_{\beta\alpha}(y)}
\]
on $U_\alpha\cap U_\beta\neq\emptyset$.
Furthermore, $h_\alpha(x,x)=h_\alpha(x)$,
$\bar\pa_x h_\alpha(x,y)$ and $\pa_y h_\alpha(x,y)$
vanish at $x=y$ to infinite order.

Let $x,y\in L^*|_{U_\alpha}$ whose local coordinates
are
$(z,v)$ and $(w,v')$, respectively. Define a
global function
\[
\psi(x,y)=
\psi(z,v,w,v')=\frac 1i(\frac{1}{h_\alpha(z,w)}v\bar
v'-1).
\]
If $x,y\in X$, write
\[
v=\sqrt{h(z)}e^{i\theta},\qquad 
v'=\sqrt{h(w)}e^{i\theta'},
\]
where $\theta,\theta'$ are real numbers.
Thus on $X$, we have
\begin{equation}\label{2-1}
\psi(x,y)=
\psi(z,\theta,w,\theta')=\frac 1i\left(
\frac{\sqrt{h(z)}\sqrt{h(w)}}{h(z,w)}e^{i(\theta-\theta')}-1\right).
\end{equation}

Let $\hat S_1,\cdots,\hat S_d$ be functions on $H^2(X)$
such that
\[
\frac{\pa}{\pa\theta}\hat S_i=\sqrt{-1}m\hat S_i
\]
for $1\leq i\leq d$. 
Using ~\cite{sz}, we can identify the functions
$\hat S_i$ to the holomorphic sections $S_i$ of $L^m$.
We assume that 
$\hat S_1,\cdots,\hat S_d$ is an orthonormal set with
respect to the measure
$d\mu=\frac{1}{n!}\pi^*\omega_g^n\wedge d\theta$. That
is
\[
\int_X\hat S_i\bar{\hat S_j} d\mu=\delta_{ij}
\]
for $1\leq i,j\leq d$. Define
\begin{equation}\label{pix}
\Pi_m(x,y)={2\pi}\sum_{i=1}^N \hat
S_i(x)\bar{\hat S_j(y)}.
\end{equation}
The basic identity related $\Pi_m$ and $\Pi$
is the following (cf. ~\cite{sz})
\begin{equation}\label{2-2}
\Pi_m(x,y)=\int_{S^1}\Pi (x,r_\theta
y)e^{\sqrt{-1}m\theta}d\theta.
\end{equation}
It is also well known
 (cf. ~\cite{BFG}~\footnote{In their paper, the
result was written in terms of   strongly
pseudoconvex domains of $\C^n$. But it is also true for
strongly pseudoconvex manifolds.}) that for the
pseudoconvex manifold
$X$, 
the Szeg\"o kernel can be written as
\begin{equation}\label{2-3}
\Pi(x,y)=\frac{u(x,y)}{\psi(x,y)^{n+1}}
+v(x,y)\log\psi(x,y), 
\end{equation}
where $u, v$ are $C^\infty$ functions defined on
$D\times D$.

Using the above notations, Zelditch~\cite{sz} (and
Catlin~\cite{Cat} for the similar result for the
Bergman kernel)  proved the following:

\begin{theorem}\label{zelditch}
With the notations as above, we have the following asymptotic
expansion:
\begin{equation}\label{sz1}
\Pi_m(x,x)\sim
m^n(a_0(x)+\frac{a_1(x)}{m}+\frac{a_2(x)}{m^2}+\cdots),
\end{equation}
where $a_i(i\geq 1)$ are smooth functions on $M$ and 
$a_0(x)=1$. The expansion is convergent in the sense
that
\begin{equation}\label{sz2}
||\Pi_m(x,x)-m^n(1+\frac{a_1(x)}{m}+\cdots+
\frac{a_k(x)}{m^k})
||_{C^s}\leq C\frac{1}{m^{k+1}},
\end{equation}
where $C$ is a constant depending on $k$, $l$ and the
manifold $M$. 
\end{theorem}

\qed

Using his result, we  prove the
following~\cite{Lu10}
\begin{theorem}\label{thm22}
The coefficients $a_i$ can be written as polynomials
of the curvature and their derivatives of $M$.
The Weyl weight of $a_k$ is $2k$ for
$k=1,2,\cdots$. In particular, we have
\[
\left\{
\begin{array}{l}
a_0=1\\
a_1=\frac 12 \rho\\
a_2=\frac 13\Delta\rho+\frac{1}{24}(|R|^2-4|Ric|^2+3\rho^2)\\
a_3=\frac 18\Delta\Delta\rho+\frac 1{24}\,div\,div\, (R,Ric)
-\frac 16 div\,div (\rho Ric)\\
+\frac{1}{48}\Delta (|R|^2-4|Ric|^2+8\rho^2)
+\frac{1}{48}\rho(\rho^2-4|Ric|^2+|R|^2)\\
+\frac{1}{24}(\sigma_3(Ric)-Ric(R,R)-R(Ric,Ric)),
\end{array}
\right.
\]
where $R, Ric$ and $\rho$ represent the curvature tensor, the Ricci
curvature and the scalar curvature of $g$, respectively
and $\Delta$ represents the Laplacian of $M$.
For the precise definition of the terms in the
expression of $a_3$, see  ~\cite[Section 5]{Lu10}.
\end{theorem}

For the above settings,
the famous Ramadanov
Conjecture~\cite{rama} (in terms of the Szeg\"o
kernel)  states that if the function $v(x,y)$ in
~\eqref{2-3} is identically zero,
then the manifold $X$ must be the sphere.

The main result of this section is the following

\begin{theorem}\label{thm23}
Let $X$ be the unit circle bundle of $L^*$ over $M$. 
If $v(x,y)$ (i.e., the log
term) of the Szeg\"o kernel of $X$ vanishes, then 
the coefficients $a_k$ in Theorem~\ref{zelditch}
vanish
for $k>n$.
\end{theorem}

{\bf Proof.} 
Using ~\eqref{2-1} and~\eqref{2-2}, if $v\equiv 0$
in~\eqref{2-3}, then we have the following identity
\[
\Pi_m(x,x)=\int_{S^1}\frac{({\sqrt{-1}})^{n+1}u(x,r_\theta
x)}{(e^{-\sqrt{-1}\theta}-1)^{n+1}}
e^{\sqrt{-1}m\theta} d\theta.
\]
We shall prove that the above expression expands to a
polynomial of $m$. Let $b>1$ be a real number, then
the above integration is understood as
\[
\Pi_m(x,x)=\underset{b\rightarrow 1}{\lim}
\int_{S^1}\frac{(\sqrt{-1})^{n+1}u(x, r_\theta x)}
{(e^{-\sqrt{-1}\theta}-b)^{n+1}}e^{\sqrt{-1}m\theta}
d\theta.
\]
Using the integration by parts $n$ times, we get
\[
\Pi_m(x,x)=
\underset{b\rightarrow 1}{\lim}
\int_{S^1}\frac{\xi(x,\theta,m)}
{e^{-\sqrt{-1}\theta}-b}e^{\sqrt{-1}m\theta}
d\theta,
\]
where $\xi(x,\theta,m)$ is a polynomial of $m$ and the
coefficients are smooth functions of $x$ and $\theta$.
By the Riemann-Lebesgue lemma, we know that that above
expression has the same asymptotic expansion as 
\[
\Pi_m(x,x)=\xi(x,0,m)\,\cdot
\underset{b\rightarrow 1}{\lim}
\int_{S^1}\frac{1}
{e^{-\sqrt{-1}\theta}-b}e^{\sqrt{-1}m\theta}
d\theta.
\]
Thus there is a {\sl polynomial} $P(x,m)$ of $m$ of
degree less than or equal to $m$ such that
\[
\Pi_m(x,x)\sim P(x,m)
\]
in the sense that 
\[
|\Pi_m(x,x)-P(x,m)|<\frac{C}{m^k}
\]
for any $k$.
Comparing the above result with the expansion in
Theorem~\ref{zelditch}, we get the conclusion.

\qed

We can also prove the above result using the stationary
phase theorem as follows:
the theorem of Boutet de Monvel and
Sj\"ostrand~\cite[Theorem 1.5 and \S 2.c]{BS1} (see
also~\cite{sz}) states that
there exists
a
symbol $s\in S^n(X\times X\times\R^+)$ of the type
\[
s(x,y,t)\sim\sum_{k=0}^\infty t^{n-k}s_k(x,y)
\]
so that
\begin{equation}\label{B-S}
\Pi(x,y)=\int_0^\infty e^{it\psi(x,y)} s(x,y,t) dt,
\end{equation}

Using~\eqref{2-2}, we have~\cite{sz}
\begin{equation}\label{sz}
\Pi_m(x,x)\sim\sum_{k=0}^\infty
m^{n-k+1}\int_0^\infty\int_{S^1} e^{im(\frac
ti(e^{i\theta}-1) -\theta)} t^{n-k} s_k(r_\theta x,x)\,
dtd\theta
\end{equation}

By the stationary phase method~\cite[Theorem
7.7.5]{Hor1}, we have 
\begin{align*}
&\int_0^\infty\int_{S^1} e^{im(\frac ti(e^{i\theta}-1)-\theta)}t^{n-k} 
s_k(r_\theta x,x) dtd\theta\\
&\sim 2\pi\sum_j m^{-j-1}L_j(t^{n-k} s_k(r_\theta x,x)),
\end{align*} 
where
\begin{align*}
&L_j(t^{n-k} s_k(r_\theta x,x))\\
&=\sum_{\nu-\mu=j}\sum_{2\nu\geq 3\mu}
i^{-j}2^{-\nu}\left(2\frac{\pa^2}{\pa t\pa\theta}
-i\frac{\pa^2}{\pa t^2}\right)^\nu
(g^\mu t^{n-k}s_k(r_\theta x,x))/{\mu!\nu!}
\end{align*}
for
\[
g(t,\theta)=\frac ti(e^{i\theta}-1)-\theta-2(t-1)\theta-i\theta^2.
\]
(cf. ~\cite{sz}). 
If $j>n-k$, then
\[
(\frac{\pa}{\pa t})^\nu (g^\mu t^{n-k} 
s_k(r_\theta x,x))\equiv 0.
\]
Thus if $j>n-k$, $L_j(t^{n-k}s_k(r_\theta
x,x))\equiv 0$. From~\eqref{sz}, we see that
\[
\Pi_m(x,x)\sim C\sum_{j+k\leq n} m^{n-k-j} L_j(t^{n-k} s_k(r_\theta
x,x)).
\]
Comparing to ~\eqref{sz1}, we have
 $a_k=0$ for $k>n$.

\qed

For the Bergman kernel of
$D$, we have the parallel result:

\begin{theorem}
If the log term of the Bergman kernel $B(x,y)$
of $D$ vanishes, then the coefficients $a_k$ in 
Theorem~\ref{zelditch} vanish for $k>n$.
\end{theorem}

\qed

\section{Order of the coefficients}

%2003.02.18.9.28

Let $d=d_m=\dim_{\C} H^0(M,L^m)$ for a fixed
integer $m$. Let $\{S_0,\cdots, S_{d-1}\}$
be
a basis of $H^0(M,L^m)$. The metrics $(h,\omega_g)$ 
define the $L^2$
inner product $(\,,\,)$ on $H^0(M,L^m)$ as
\[
(S_A,S_B)=\int_M<S_A,S_B>dV_g,
\quad A,B=0,\cdots,d-1,
\]
where $<S_A,S_B>$ is the pointwise inner product with respect
to $h_m$ and $dV_g=\frac{1}{n!}\omega_g^n$.
If $T_0,\cdots, T_{d-1}$ is an orthonormal basis
of $H^0(M,L^m)$
with respect to the above inner product, then
 we define the sum of the pointwise norm
\[
||T_0||^2(z)+\cdots+||T_{d-1}||^2(z)
=<T_0,T_0>(z)+\cdots +<T_{d-1}, T_{d-1}>(z)
\]
to be the Bergman potential of the metric. The key link
between the Bergman potential and the Szeg\"o kernel
is (cf.~\cite{sz})
\[
\Pi_m(x,x)=||T_0||^2(z)+\cdots+||T_{d-1}||^2(z),
\]
where $\pi(x)=z$, and $\Pi_m(x,x)$ is defined
in~\eqref{pix}.

Now we assume $S_0, \cdots, S_{d-1}$ is a basis of the
space $H^0(M,L^m)$.
We further assume that at a  point $z \in M$,
\[
S_0(z)\neq 0,\quad S_{A}(z)=0,\quad A=1,\cdots, d-1.
\]
Suppose 
\[
F_{AB}=(S_A,S_B), \quad A, B=0,\cdots, d-1.
\]
Then $(F_{AB})$ is the metric matrix which is positive
Hermitian. Let
$(I_{AB})$ be the inverse matrix of $(F_{AB})$.
Let $x\in X$ such that $\pi(x)=z$. Then by linear algebra we have
(cf. ~\cite{Lu10})
\begin{equation}\label{import}
\Pi_m(x,x)=I_{00} ||S_0(z)||^2_{h_m},
\end{equation}
where $||S_0(z)||_{h_m}^2=<S_0(z), S_0(z)>$ is the pointwise
norm of the section $S_0(z)$. 

The main result of this section is to prove that all the coefficients
$a_k$ in the theorem of Zelditch
(Theorem~\ref{zelditch}) can be represented by
$C\Delta^{k-1}\rho$ plus ``lower order terms", where
$C\neq 0$ is a constant depending only on $k$ and $n$ and $\rho$ is the
scalar curvature of $M$. To make the above
statement rigorous,
we need the following definition:

\begin{definition}\label{md1}
Let $R$ be a component of the $i$-th order covariant derivative of the
curvature
tensor, or the Ricci tensor, or the scalar curvature
at a fixed point
where $i\geq 0$. Define the weight
$w(R)$ and the order $ord(R)$ of $R$ to  be the number $(1+\frac i2)$
and $\frac i2$, respectively. For example,
\[
\begin{array}{l}
w(\crr ijkl)=w(\css ij)=w(\rho)=1, \\
ord(\crr ijkl)=ord(\css ij)=ord(\rho)=0
\end{array}
\]
and
\[
\begin{array}{l}
w(R_{i\bar jk\bar l,m})=\frac 32,\\
ord(R_{i\bar jk\bar l,m})=\frac 12.
\end{array}
\]
In particular, the weight and the order of a constant are zero.
The concepts of weight and order can be extended to  monomials of the
curvature and its derivatives by assuming that
\begin{align*}
& w(f_1f_2)=w(f_1)+w(f_2),\\
& ord(f_1f_2)=ord(f_1)+ord(f_2),
\end{align*}
where $f_1,f_2$ are monomials.
If $f=\sum f_i$ with $f_i$ monomials of the same
weight or order, then we define $w(f)=w(f_1)$ and
$ord(f)=ord(f_1)$, respectively.
\end{definition}

\begin{rem} The definition of the weight here is
half of the Weyl weight in Fefferman's paper~\cite
{F1}.
\end{rem}

Let $A$ be the set of all monomials of the curvature and its derivatives at
a fixed point $z\in M$. Define
\[
A'=\{f\in A| ord(f)\leq w(f)-2\}.
\]
Let $B$ and $B'$ be the complex vector spaces generated by $A$ and $A'$,
respectively.
We have the following simple relation between  weight and order:

\begin{lemma}\label{sp}
 For any $f_1, f_2\in A$ with $w(f_1),w(f_2)\neq 0$, we have
$f_1f_2\in A'$. In particular, $B'$ is an ideal of $B$.
\end{lemma}

{\bf Proof.} If  $w(f_1),w(f_2)\neq 0$, then we have
\[
ord(f_i)\leq w(f_i)-1,\qquad i=1,2.
\]
Thus
\[
ord(f_1f_2)=ord(f_1)+ord(f_2)\leq w(f_1)+w(f_2)-2=w(f_1f_2)-2.
\]

\qed

The main result of this section is

\begin{theorem}\label{thm31}
With the notations as above, for any $k\geq 1$, there is a constant
$C=C(k,n)\neq 0$ such that
\[
a_k\equiv C\Delta^{k-1}\rho\qquad mod\,(B'),
\]
where $\rho$ is the scalar curvature of $M$ and $\Delta$ is the Laplace
operator of $M$.
\end{theorem}

In order to prove the theorem, we must estimate the quantities in
~\eqref{import}. We
construct peak sections of $L^m$ for $m$ large. So let's quickly review
the concept of peak sections which was initiated in~\cite{T3}.

Choose a local normal
coordinate $(z_1,\cdots,z_n)$ centered at $z$  such that
the Hermitian matrix $(g_{\alpha\bar\beta})$ satisfies
\begin{align}
\begin{split}
& g_{\alpha\bar\beta}(z)=\delta_{\alpha\beta},\\
&\frac{\pa^{p_1+\cdots+p_n}g_{\alpha\bar\beta}}{\pa z_1^{p_1}\cdots\pa
z_n^{p_n}}
(z)=0
\end{split}
\end{align}
for $\alpha, \beta=1,\cdots,n$ and any nonnegative integers $p_1,\cdots,
p_n$ with $p_1+\cdots+p_n\neq 0$. 
\footnote{Note that there is a little ambiguity about the notation here.
We use $z$ to denote both the point and its local coordinates. However,
it should be clear from the context.}
Such a local
coordinate
system, which is known as the $K$-coordinate system
(cf. 
~\cite{Bo} or~\cite{ru} for
details) exists and is unique up to an affine
transformation. We choose a local holomorphic frame
$e_L$ of $L$ at $z$ such that the local representation
function $a$ of the Hermitian metric
$h$ has the properties
\begin{equation}
a(z)=1, \frac{\pa^{p_1+\cdots+p_n}a}{\pa z_1^{p_1}\cdots\pa
z_n^{p_n}}(z)=0
\end{equation}
for any nonnegative integers $(p_1,\cdots,p_n)$ with $p_1+\cdots+p_n\neq
0$.

Suppose that the local coordinate $(z_1,\cdots,z_n)$ is defined on an
open neighborhood $U$ of $x_0$ in $M$. Define
the function $|z|$
by
$|z|=\sqrt{|z_1|^2+\cdots+|z_n|^2}$  for $z\in U$.

Let $\Z_+^n$ be the set of $n$-tuple of integers $(p_1,\cdots,p_n)$
such that $p_i\geq 0 (i=1,\cdots, n)$. Let $P=(p_1,\cdots,p_n)$.
Define
\begin{equation}\label{zp}
z^P=z_1^{p_1}\cdots z_n^{p_n}
\end{equation}
and
\[
p=|P|=p_1+\cdots+p_n.
\] 

The following lemma
is  proved in ~\cite{T5}
using
the standard $\bar\pa$-estimates (see e.g.~\cite{Hor}).

\begin{lemma}\label{tian}
For  $P=(p_1,\cdots,p_n)\in \Z_+^n$, and an
integer $p'>p=p_1+\cdots+p_n$, there exists an $m_0>0$ such that for
$m>m_0$, there is a holomorphic global section
$S^{p'}_{P,m}$ in $H^0(M,L^m)$, satisfying
\[
\int_M||S^{p'}_{P,m}||^2_{h_m} dV_g=1,\qquad
\int_{M\backslash\{\rr\}}
||S^{p'}_{P,m}||^2_{h_m} dV_g=\oo{2p'},
\] 
and $S^{p'}_{P,m}$ can be decomposed as
\[
S^{p'}_{P,m}=\tilde S_{P,m}+u_{P,m},\qquad (\tilde
S_{P,m}
\,
\text{and}\,u_{P,m}\,
\text{not necessarily continuous})
\]
such that
\[
\tilde S_{P,m}(x)=\left\{
\begin{array}{ll}
\lambda_Pz^Pe_L^m(1+\oo{2p'}) &x\in\{ \rr\}\\
0                        &x\in M\backslash\{\rr\},
\end{array}
\right.
\]
\[
u_{P,m}(x)=\ooo{2p'}\qquad x\in U,
\]
and
\[
\int_M||u_{P,m}||_{h_m}^2 dV_g=\oo{2p'},
\]   
where
$\oo{2p'}$ denotes a quantity dominated by $C/m^{2p'}$ with the constant
$C$ depending only on $p'$ and the geometry of $M$. Moreover
\[
\lambda_P^{-2}   
=\int_\rr |z^P|^2a^mdV_g.
\]
\end{lemma}

\qed

Define an order $\geq$ on the multiple indices $P$
as follows:  $P\geq Q$, if
\begin{enumerate}
\item $|P|>|Q|$ or,
\item $|P|=|Q|$ and $p_j=q_j$ but $p_{j+1}>q_{j+1}$ for some $0\leq j
\leq n$.   
\end{enumerate}
Using this order, there is a one-one order preserving correspondence
$\kappa$
between $\{0,1,2,\cdots\}$ and $\{P| P\in \Z^n_+\}$.

We need the following proposition in~\cite{Lu10}:

\begin{prop}\label{spq}
We have the following expansion for any
$p'>t+2(n+p+q)$,
\[
(S_{P,m}^{p'},S_{Q,m}^{p'})
=\frac{1}{m^\delta}(a_0+\frac{a_1}{m}+\cdots+\frac{a_{t-1}}{m^{t-1}}+\oo{t}),
\]
where $\delta=1$ or $1/2$ and
where all the $a_i$'s are polynomials of the curvature and its derivatives
such that
\[
ord (a_i)=i+\delta.
\]
\end{prop}

\qed

It has been proved in~\cite[Theorem 3.1]{Lu10} that for any $t>0$, there
is an
$s>0$ such that up to 
$\oo{t}$, $I_{00}$ depends only on $S_{\kappa(i)}\,(i=0,\cdots,s)$.
More precisely, let
\[
F_{ab}'=(S_{\kappa(a)},S_{\kappa(b)}),\quad
0\leq a,b\leq s.
\]
Let $(I'_{ab})$ be the inverse matrix of $(F_{ab}')$, then
\[
I_{00}=I_{00}'+\oo{t}.
\]
Let
\[
(F_{ab}')=
\begin{pmatrix}
1 & M_{21}\\
M_{12} & M_{22}
\end{pmatrix},
\]
where $M_{12}\in \C^s$, $M_{21}^T\in \C^s$ and
$M_{22}$ is an $s\times s$ matrix. By an elementary computation we have
\begin{equation}\label{ele}
I_{00}'=1+M_{12}^T(M_{22}-M_{21}M_{12})^{-1}
M_{21}.
\end{equation}

In~\cite[Lemma 2.2]{Lu10}, we know that $M_{12}=\oo{}$. In particular,
for any monomial in any entry $e$ of $M_{12}$,
$ord(e)\leq w(e)-1$. Thus by ~\eqref{ele},
Proposition~\ref{spq} and Lemma~\ref{sp} we have
\begin{equation}\label{3-6}
I_{00}\equiv 1\qquad mod\,(B').
\end{equation}

\smallskip

Next let's consider $||S_0(z)||_{h_m}^2$. By
Lemma~\ref{tian} we see that
\[
||S_0(z)||^2_{h_m}=\lambda_0^{2}+O(\frac{1}{m^N})
\]
for any $N$, 
where
\[
\lambda_0^{-2}=\int_\rr a^m dV_g.
\]
Let $\xi=\log a+|z|^2$, $\eta=\log\det
g_{\alpha\bar\beta}$, we have
\[
\lambda_0^{-2}=\int_\rr e^{m\xi+\eta} 
e^{-m|z|^2} dV_0,
\]
for the Euclidean volume form $dV_0$.

By Lemma~\ref{sp}, we see that
\[
\lambda_0^{-2}\equiv\int_\rr (1+m\xi+\eta) e^{-m|z|^2} dV_0
\qquad mod\,(B').
\]

Using the fact
\[
\int_{\C^n}|z^{p_1}\cdots z^{p_n}|^2e^{-m|z|^2}dV_0
=\frac{p_1!\cdots p_n!}
{m^{n+p}},
\]
we have the following
\begin{align}\label{3-1}
\begin{split}
&\lambda_0^{-2}\equiv
\frac{1}{m^n}\\
&+\frac{1}{m^n}
\sum_{k=1}^N( \frac{1}{(k+1)!}
\Delta_c^{k+1}\xi+\frac{1}{k!}
\Delta_c^j\eta)\frac{1}{m^k}+\oo{N+n}\qquad mod\,(B')
\end{split}
\end{align}
for any $N$, where $\Delta_c$ is the complex Laplace operator
 on $\C^n$, defined by
\begin{equation}\label{3-1-1}
\Delta_c=\sum_{i=1}^n\frac{\pa^2}{\pa z_i\pa\bar z_i}.
\end{equation}
As before, $\Delta$ will be the Laplace operator on
$M$. It is not hard to see that
\begin{equation}\label{3-2}
\Delta_c^j\eta\equiv-\Delta^{j-1}\rho\qquad mod\,(B').
\end{equation}
for $j\geq 1$. Using the same method, we have
\begin{equation}\label{3-3}
\Delta_c^{j+1}\xi\equiv \Delta^{j-1}
\rho\qquad mod\,(B').
\end{equation}
Combining~\eqref{3-1},~\eqref{3-2} and~\eqref{3-3},
we have
\[
\lambda_0^{-2}\equiv\frac{1}{m^n}(1-\sum_{k=1}^N
\frac{k}
{(k+1)!m^k}\Delta^{k-1}\rho)+\oo{N+n}\qquad mod\,(B').
\]
Thus
\begin{equation}\label{3-5}
\lambda_0^2\equiv
m^n(1+\sum_{k=1}^N\frac{k}{(k+1)!}\frac{1}{m^k}
\Delta^{k-1}\rho)+\oo{N+n}\qquad mod\,(B').
\end{equation}
Comparing the above equation with~\eqref{3-6}, we have
\[
a_k=\frac{k}{(k+1)!}\Delta^{k-1}\rho\qquad mod
(B').
\]
and Theorem~\ref{thm31} is proved.

\qed

For the rest of this paper, we will study the
coefficients $a_k$ in the Tian-Yau-Zelditch
expansion for different metrics. We shall thus use the
notation $a_k(x,h)$, where $x\in M$ and $h$ is the
Hermitian metric on $L$, to explicitly
represent the dependence of  the
coefficients to the metric.

Let $\ph\in C^\infty(M)$. Then $he^\ph$ defines
a Hermitian metric on $L$. Let $\omega_\ph=
\bb(g_{i\bar
j}-\pa_i\bar\pa_j\ph)dz_i\wedge d\bar z_j$ be the
corresponding \ka form. We have the following

\begin{cor}\label{cor31}
Using the notation as in  Theorem~\ref{thm31}, let 
\[
\omega_\ph>\frac 12
\omega.
\]
Then if
\[
a_k(x,he^\ph)=0,
\]
then there is a constant $C(l)$, depending on the 
$C^{l-4}$ bound of the curvature of $\omega$, such
that
\[
||\ph||_{C^l}\leq C(l)||\ph||_{C^{2k+2}}
\]
for $l>2k+2$.
\end{cor}

{\bf Proof.} We have
\[
a_k\equiv C\Delta^{k-1}\rho\quad ({\rm mod} B')
\]
where
$\rho$ is the scalar  curvature of $\omega_\ph$.
By the Schauder estimate we have
\[
||\rho||_{C^{2k-2,\frac 12}}\leq
C||\phi||_{C^{2k+1,\frac 12}}.
\]
Since 
\[
\rho=-\Delta\log\frac{\omega_\phi^n}{\omega^n}
-g_\phi^{i\bar j}\pa_i\bar\pa_j\log\omega^n,
\]
we have
\[
||\frac{\omega_\phi^n}{\omega^n}||_{C^{2k,\frac
12}}\leq C||\phi||_{C^{2k+1,\frac 12}}.
\]
By the Schauder estimate again we have
\[
||\ph||_{C^{2k+2,\frac 12}}\leq C||\ph||_{C^{2k+1,\frac
12}}.
\]
The boot strapping method gives the higher order
estimates.

\qed

\section{The uniformity of the expansion}

% 2002.2.18.9.42

As in the previous sections, 
let $h$ be a Hermitian metric on the line bundle $L$
over $M$. Let
$\ph$ be a smooth function such that $h_t=he^{t\ph}$
be a family of Hermitian metrics with $-\pa\bar\pa\log
h_1>0$. Assume
that
$S_0,\cdots,S_{d-1}$ is a  basis for the Hermitian
vector space $H^0(M,L^m)$ such that at a fixed point $x$, 
$S_0(x)\neq 0$ but $S_j(x)=0$ for $j\neq 0$. We use
$<S_i,S_j>_t$ and $(S_i, S_j)_t$ to denote the pointwise and
the $L^2$ inner product respectively with respect to the
metric $h_t$. Let $F_{t,\alpha\beta}=(S_\alpha, S_\beta)_t$.
Let $I_{t,\alpha\beta}$ be the  inverse
matrix of 
$F_{t,\alpha\beta}$. Then by~\eqref{import}, the Bergman 
potential at $x$ with respect to the metric $h_t$ is
\begin{equation}\label{5-1}
\sigma(t)=I_{t,00}||S_0||^2_t(x).
\end{equation}

The following result was pointed out
by Zelditch (cf. ~\cite[Proposition 6]{SKD1}).  

\begin{prop}
In Theorem~\ref{zelditch},  for any $s$ and $k$, there is
a number $N=N(s,k)$ such that if the metric $\omega$ is
lowerly bounded  and is bounded in $C^N$ by a constant $C_1$
with some reference metric,
then the constant $C$ in ~\eqref{sz2} depends only on $s,k$
and the constant $C_1$.
\end{prop}

\qed

Let  $\ph\in C^\infty$. Then in the expansion of 
$\sigma(t)$, the constant $C$ is independent
of $t$ for $0\leq t\leq \frac 12$ by the above proposition.

It is natural to ask whether
one can take the derivative of the
expansion of $\sigma(t)$ in order to get the expansion of 
$\sigma'(0)$. The main result of this section is to confirm
that this is indeed the case.

Let $I_{\alpha\beta}=I_{0,\alpha\beta}$. 
We have the  following
\begin{prop}\label{pp51}
At the fixed point $x$, using the above notations, we have
\begin{equation}\label{prop51}
\sigma'(0)=-\sigma(0)
I_{00}^{-1}\int_M
(m\ph-\Delta\ph)I_{o\alpha}<S_\alpha,S_\beta>
I_{\beta 0}\,\frac{\omega_0^n}{n!},
\end{equation}
where we assume that
$\ph(x)=0$  without losing generality. 
\end{prop}

{\bf Proof.} 
For fixed $m$,
we have
\[
F_{t,\alpha\beta}=(S_\alpha, S_\beta)_t
=F_{0,\alpha\beta}+t\int_M (m\ph-\Delta\ph)<S_\alpha,
S_\beta>\frac{\omega_0^n}{n!}+O(t^2).
\]

A straightforward computation gives
\[
I_{t,00}=I_{00}-t\int_M(m\ph-\Delta\ph)I_{0\alpha}
<S_\alpha, S_\beta>I_{\beta 0}\,\frac{\omega_0^n}{n!}+O(t^2).
\]
We also have
\[
||S_0||^2_t=||S_0||^2.
\]
The lemma follows from ~\eqref{5-1} and the above two
equations.

\qed

In order to get the uniform estimate, we have to establish a
uniform version of the above proposition. The main difficulty
here is
that the size of the matrix is very  large (of the size
$m^n$). The technique we use here is to  choose a special kind
of basis under which the matrix $F_{\alpha\beta}$ takes the
form of  ~\eqref{uv}.

We make the following definition from (cf. ~\cite{Lu10}):

\begin{definition}
We say $N=\{N(m)\}$ is a sequence of $s\times s$ block
matrices with block
number $t\in \Z$, if for each $m$,
\[
N=N(m)=
\begin{pmatrix}
N_{11}(m)  & \cdots & N_{1t}(m)\\
\vdots  &\ddots  & \vdots\\
N_{t1}(m) &\cdots & N_{tt}(m)
\end{pmatrix}
\]
such that for $1\leq i,j\leq t$, $N_{ij}$ is a
$\sigma(i)\times\sigma(j)$
matrix and
\[
\sum_{i=1}^t\sigma(i)=s,
\]
where $\sigma:\{1,\cdots, t\}\rightarrow \Z_+$ assigns each number
in $\{1,\cdots,t\}$ a positive integer. We say that
$\{N(m)\}$ is of type $A(p)$ for a positive integer $p$, if
for any entry
$s$ of the matrix $N$, we have
\begin{enumerate}
\item
If $s$ is a diagonal entry of $N_{ii} (1\leq i\leq t)$ ,
then we have the the following
Taylor expansion
\[
s=1+\frac{s_1}{m}+\cdots+\frac{s_{p-1}}{m^{p-1}}+\oo{p}.
\]
\item
If $s$ is not a diagonal entry of $N_{ii} (1\leq i\leq t)$,
then we have the Taylor expansion
\[
s=\frac{1}{m^\delta}(s_0+\frac{s_1}{m}+\cdots+\frac{s_{p-1}}{m^{p-1}}
+\oo{p}),
\]
where $\delta$ is equal to 1 or $\frac 32$. 
\item
If $s$ is an entry of the matrix $N_{ij}$ for which
$|i-j|=1$, then
$s=\oo{\frac 32}$. In addition, if $i\neq t$ or $j\neq t$, then
we have the Taylor expansion
\[
s=\frac{1}{m^\delta}(s_0+\frac{s_1}{m}+\cdots+\frac{s_{p-1}}{m^{p-1}}
+\oo{p}),
\]
where $\delta$ is equal to 1 or $\frac 32$.
\item
If $s$ is an entry of $N_{ij}$ for which $|i-j|>1$, then
\[
s=\oo{p}.
\]
\end{enumerate}
The set of all quantities $(s_1/m,\cdots,s_{p-1}/m^{p-1})$,
or $(\frac{s_0}{m^\delta},
\frac{s_1}{m^{1+\delta}},\cdots, \frac{s_{p-1}}{m^{p+\delta-1}})$ for
 $s$ running from all the entries 
of $N_{ij}$ where $|i-j|\leq 1$ and $i\neq t$ or $j\neq t$ 
are called
the Taylor Data of order $p$.
\end{definition}

\begin{rem}
Since $N_{ij}=\oo{p}$ for $|i-j|>1$, it can be treated as zero
when we are only interested in the expansion of order up to
$p-1$, and when the rank of the matrix is bounded by a
constant depending only on $p$. A matrix whose entries
$N_{ij}=0$ for
$|i-j|>1$ is called a tri-diagonal matrix. For such a matrix,
we have a simple iteration process for finding its inverse
matrix~\cite{Lu10}.
\end{rem}

\begin{prop}\label{pp52}
For any positive integer $p$, there is a number $\xi(p)$
such that there is a $\xi(p)\times \xi(p)$ block matrix
$N$ of block number $(p+1)$. $N$ is of type $A(p)$.
Furthermore, the matrix $(F_{AB})$ can be represented as
\[
(F_{AB})=\left(
\begin{array}{cc}
N& 0\\
0 &E(d-\xi(p))
\end{array}
\right),
\]
where $E(d-\xi(p))$ is the $(d-\xi(p))\times (d-\xi(p))$
identity matrix.
\end{prop}

{\bf Proof.}
We construct such a matrix using the peak sections in
Lemma~\ref{tian}.  For a multiple indices, define
$|P|=p_1+\cdots+p_n$. Suppose 
\[
V_k=\{S\in H^0(M,L^m)|D^QS(x_0)=0\,\text{for}\, |Q|\leq k\}
\]
for $k=1,2,\cdots$, where $Q\in \Z_+^n$ is a multiple indices, and
$D$ is a covariant derivative on the bundle $L^m$. 
$V_k=\{0\}$ for $k$ sufficiently large.
For fixed $p$, let $p'=n+8p(p-1)$. Suppose that $m$ is large enough such
that $H^0(M,L^m)$ is spanned by the $S_{P,m}^{p'}$'s for the
multiple indices $|P|\leq 2p(p-1)$ and $V_{2p(p-1)}$. Let $r=d-\dim
V_{2p(p-1)}$.
Then $r$ only depends on $p$ and $n$.
Let 
$T_1,\cdots,T_{d-r}$ be an orthonormal basis
of $V_{2p(p-1)}$
such that 
\[
(S_{P,m}^{p'},T_\alpha)=0
\]
for $|P|\leq 2p(p-1)$ and $\alpha>r$. Let $s(k)=\dim V_k$
for $k\in \Z$. 
For any $1\leq i,j \leq p$, let
$N_{ij}$ be the matrix formed by
$(S_{P,m}^{p'},S_{Q,m}^{p'})$ where
$2p(i-2)\leq|P|\leq 2p(i-1)$ and $2p(j-2)\leq|Q|\leq 2p(j-1)$.
Furthermore,
define $N_{i(p+1)}$ to be the matrix whose entries are
$(S_P,T_\alpha)$ for $2p(i-2)\leq|P|\leq 2p(i-1)$ and
$1\leq\alpha\leq r$. Define
$N_{(p+1)i}$ to be the complex conjugate of  $N_{i(p+1)}$.
Finally, define
$N_{(p+1)(p+1)}$ to be the $r\times r$ unit matrix $E(r)$.
Then it is easy to check that $N=(N_{ij})$ is 
a sequence of block matrices of type $A(p)$ with
the block number $p+1$ by using the result of Ruan (cf.
~\cite[Lemma 2.2]{Lu10} and Proposition~\ref{spq}.
Let $\xi(p)=2r=2(d-\dim V_{2p(p-1)})$. Then $\xi(p)$, which is
the rank of the matrix $N$, depends only on $p$ and $n$.

Define an order $\geq$ on the multiple indices $P$ 
as follows:  $P\geq Q$, if
\begin{enumerate}
\item $|P|>|Q|$ or;
\item $|P|=|Q|$ and $p_j=q_j$ but $p_{j+1}>q_{j+1}$ for some $0\leq j
\leq n$.
\end{enumerate}
Using this order, there is a one-one order preserving correspondence $\kappa$
between $\{0,\cdots, r-1\}$ and $\{P||P|\leq 2p(p-1)\}$.

Define
\[
S_A=\left\{
\begin{array}{ll}
S^{p'}_{\kappa(A),m} & A\leq r-1\\
T_{A-r+1}  & A\geq r
\end{array}
\right..
\]

Comparing the matrix $N$ to the metric matrix $F_{AB}=((S_A,
S_B)), (A,B=0,\cdots, d-1)$, by the choice of the basis, we
see that
\begin{equation}\label{uv}
(F_{AB})=
\begin{pmatrix}
N & 0\\
0 & E(d-2r)
\end{pmatrix},
\end{equation}
where $E(d-2r)$ is the $(d-2r)\times (d-2r)$ identity matrix.

\qed

Using the above result, we have

\begin{theorem}\label{thm52}
There is an expansion of $\sigma'(0)$
\[
\sigma'(0)\sim
m^n(b_0+\frac{b_1}{m}+\frac{b_2}{m^2}+\cdots),
\]
in the sense that for any $k$,
\[
||\sigma'(0)-m^n(b_0+\cdots+\frac{b_k}{m^k})||_{C^0}
\leq
\frac{C}{m^{k+1}},
\]
where the constant $C$ depends on $k$ and the manifold 
$M$ but is
independent to $m$.\footnote{The expansion is convergent even
in the $C^\infty$ norm, though we don't need the fact. One
may also prove the theorem using the paramatrix of the 
Szeg\"o kernel, similar to what Zelditch did in~\cite{sz}. We
may
 have to cope with the quantity
on different circle bundles if using the Szeg\"o
kernel method.}
\end{theorem}

{\bf Proof.} By~\cite{sz}, we know that there is
an asymptotic expansion of $\sigma(0)$. By~\cite[Theorem
3.1]{Lu10}, we have the asymptotic expansion of $I_{00}$.
Thus in order to give the expansion of $\sigma'(0)$, in
terms of ~\eqref{prop51}, we 
just need to prove that for any smooth function $\psi$, there
is an asymptotic expansion of the expression 
\begin{equation}\label{exp1}
\sum_{\alpha,\beta=0}^{d-1}
\int_M\psi I_{0\alpha}<S_\alpha, S_\beta> I_{\beta
0}\,\omega_0^n.
\end{equation}
We choose the basis $S_0,\cdots, S_{d-1}$ as in
Proposition~\ref{pp52}. By the proposition, we have 
$I_{0\alpha}=0$ for $\alpha>2r$, where $r$ is 
the size of the matrix $N$ in~\eqref{uv}.  For
each fixed
$\alpha,\beta$, it is easy to see that there is an
asymptotic expansion for the term
$\int_M\psi I_{0\alpha}<S_\alpha, S_\beta> I_{\beta
0}\,\omega_0^n$. The theorem
follows from the fact that  $r$ is independent of $m$.

We now prove the main result of this section:

\begin{theorem}\label{34}
Suppose that we have the following
expansion of $\sigma(0)$ for $t$ small:~\footnote{Here
$a_i(x,t)=a_i(x,he^{t\ph})$ for short.}
\[
\sigma(t)\sim m^n(a_0(x,t)+\frac{a_1(x,t)}{m}+\cdots)
\]
in the sense that
\begin{equation}\label{5-4}
||\sigma(t)-m^n(a_0(x,t)+\frac{a_1(x,t)}{m}+\cdots
+\frac{a_k(x,t)}{m^k})||_{C^0}\leq
\frac{C}{m^{k+1}},
\end{equation}
where $k\geq 1$ is an integer and $C$ is independent to $m$
and
$t$. Then the expansion  of $\sigma'(0)$, if exists, must be
of the form
\[
\sigma'(0)\sim
m^n(\left.\frac{d}{dt}\right|_{t=0}a_0(x,t)+
\left.\frac{d}{dt}\right|_{t=0}\frac{a_1(x,t)}{m}+\cdots).
\]
\end{theorem}

{\bf Proof.} 
In what follows, we denote $C$ to be a general constant
that is independent to $m$ and $t$.
Suppose that the expansion of $\sigma'(0)$ is
\[
\sigma'(0)\sim m^n(b_0(x)+\frac{b_1(x)}{m}+
\frac{b_2(x)}{m^2}+\cdots),
\]
with 
\[
||\sigma'(0)-m^n(b_0(x)+\frac{b_1(x)}{m}+\cdots+
\frac{b_k(x)}{m^k})||_{C^0}\leq\frac{C}{m^{k+1}}.
\]
Using~\eqref{5-4} and the above inequality, we have
\begin{equation}\label{4-20}
m^n||\sum_{i=1}^k\frac{1}{m^i}(\frac{a_i(x,t)-a_i(x,0)}{
t}-b_i(x))||_{C^0}\leq\frac{3C}{|t|\, m^{k+1}}+\left|
\frac{\sigma(t)-\sigma(0)}{
t}-\sigma'(0)\right|.
\end{equation}
If we choose the basis of $H^0(M,L^m)$ as in
Proposition~\ref{pp52}, then we have
\[
\left|\frac{\sigma(t)-\sigma(0)}{
t}-\sigma'(0)\right|\leq Cm^2|t|
\]
when $mt$ is small.

Thus ~\eqref{4-20} becomes
\begin{equation}\label{4-21}
m^n||\sum_{i=1}^k\frac{1}{m^i}(\frac{a_i(x,t)-a_i(x,0)}{
t}-b_i(x))||_{C^0}\leq\frac{3C}{|t|\,
m^{k+1}}+Cm^2|t|.
\end{equation}
The above 
inequality holds true for any $k, m$ and $t$ (constant
$C$ depends on $k$). If we choose $k=2$, $|t|=1/m^{(5/2)}$,
then letting $m\rightarrow\infty$, we have
\[
\left.\frac{d}{dt}\right|_{t=0}a_1(x,t)=b_1(x).
\]
Now we assume that for any $1\leq i<j$, we have
\[
\left.\frac{d}{dt}\right|_{t=0}a_i(x,t)=b_i(x).
\]
Since all $a_i(x,t)$ are all
differentiable, for $t$ small, there is a constant $C$
such that
\[
\left|\frac{a_i(x,t)-a_i(x,0)}{t}-b_i(x)\right|\leq
C|t|
\]
for $1\leq i<j$ and
\[
\left|\frac{a_i(x,t)-a_i(x,0)}{t}-b_i(x)\right|\leq C
\]
for $i\geq j$.
Assuming that $k>j$, then from ~\eqref{4-21}, we have
\[
m^{n-j}\left|\frac{a_j(x,t)-a_j(x,0)}{
t}-b_j(x)\right|\leq
\frac{3C}{m^{k+1}|t|}+Cm^2|t|+Cm^{n-1}j|t|
+Ckm^{n-j-1}.
\]
We assume $|t|=1/m^{j}$ and let $k>2j$, then the above
inequality implies the conclusion of the theorem.

\qed

\section{The general case}

%2003.2.18.9.44

In this section, we prove Theorem~\ref{fud}. First we establish
some general estimate that will be used for the rest of the paper.

We  use $a_l(x,h)$ to denote the
$l$-th coefficient in the 
Tian-Yau-Zelditch expansion (Theorem~\ref{thm22}), where
$h$ is the Hermitian metric on the bundle $L$ and $x\in M$.

\begin{lemma} 
Let $l$ be a nonnegative integer. 
Let $\omega=-\bb\pa\bar\pa\log h$.
Let $\ph\in C^{2l+2}$
satisfy
\[
\left\{
\begin{array}{l}
||\ph||_{C^{2l+2}}\leq 1,\\
\frac 12\omega+\sqrt{-1}\pa\bar\pa\ph>0.
\end{array}
\right.
\]
Then there is a constant $C$, depending on $l$ and the 
$C^{2(l-1)}$ curvature
bound of the metric $\omega$, such that
\begin{equation}\label{f-1}
|a_l(x,he^{t\ph})-a_l(x,h)-\left.\frac{d}{ds}
\right|_{s=0} a_l(x,he^{s\ph}) t|\leq C t^2
\end{equation}
for $0\leq t\leq 1$.
Furthermore, for a metric $h'$ which is $C^{2(l+1)}$
close to $h$,
we have the following inequality
\begin{equation}\label{f-2}
|a_l(x,h'e^{t\ph})-a_l(x,h')|\leq C_1t
\end{equation}
for $0\leq t\leq 1$ and for the constant $C_1$ depending
only on $l$, the $C^{2(l-1)}$ bound of the curvature of $h$,
and the $C^{2(l+1)}$ norm of $\ph$.
\end{lemma}

{\bf Proof.} By Theorem~\ref{thm22}, we know that 
$a_l(x,he^{t\ph})$ is a polynomial of Weyl weight $2l$.
That means $a_l(x,he^{t\ph})$ is a smooth function
of the curvature, its derivative of $\omega$ of degree up to
$2(l-1)$ and of $\ph$, its derivative of degree up
to $2(l+1)$. Using the assumption that $\bb
(-\frac 12\pa\bar\pa\log h+\pa\bar\pa\ph)>0$, we can 
expand $a_l(x, e^{t\ph})$ as the  Taylor  series
of $t$ with the coefficients depend on $l$, the 
$C^{2(l-1)}$ bound of the curvature of $\omega$,
and the $C^{2(l+1)}$ norm of the function $\ph$. 
~\eqref{f-1} follows from the Taylor expansion.

We note that the constant $C$ only depends on the bounds
of the curvature and the function $\ph$. 
~\eqref{f-2} follows from this observation.

\qed

{\bf Proof of Theorem~\ref{fud}.} If the theorem is not true,
then we have an infinite  dimensional vector space $V$ such that
for any
$\ph\in V$, the log term for the metric $he^{t\phi}$
is zero for $t$ small enough. By Theorem~\ref{thm23},
we have
$a_{n+1}(x,he^{t\ph})\equiv 0$.
By Theorem~\ref{thm31}, we have
\[
a_{n+1}(x,he^{t\ph})=C\Delta_t^n\rho_t+\text{lower order terms},
\]
where $\Delta_t$ and $\rho_t$ are the Laplacian and the scalar
curvature of the metric 
$\omega-t\bb\pa\bar\pa\ph$, respectively.
A straightforward computation gives
\[
\left.\frac{d}{dt}\right|_{t=0}\rho_t=\Delta^2\ph+R_{j\bar
i}\phi_{i\bar j},
\]
where $\Delta$ is the Laplacian of $\omega$.
Thus we have
\begin{equation}\label{3-8}
0=\left.\frac{d}{dt}\right|_{t=0}a_{n+1}(x,he^{t\ph})
=C\Delta^{n+2}\ph+\text{lower order terms}.
\end{equation}
Since the above identity is a linear elliptic equation of $\ph$.
The solution space is a finite dimensional space by the 
Schauder estimates.

\qed

\section{The cases of complex projective spaces}
% 2003.2.18.9.47
In this section, we study the unit circle bundle of the
universal line bundle of the complex projective space
$\cp$. First,
we prove Theorem~\ref{fud1}, which is parallel to the case of 
pseudoconvex domain in $\C^2$.

{\bf Proof of Theorem~\ref{fud2}.}
By Theorem~\ref{thm23},
we must have $a_2=0$. By Theorem~\ref{thm22}, we have
\[
a_2=\frac 13\Delta\rho+\frac{1}{24}(|R|^2-4|Ric|^2+3\rho^2).
\]
Since $n=1$, the above equation is reduced to
\[
\Delta\rho=0.
\]
Thus the scalar curvature must be constant. Since $M=\C
P^1=S^2$, the constant $\rho$ must be positive and thus
the metric must be the standard one.

\qed

We now assume that $(M,L)=(\C P^n, H)$,
where $H$ is the hyperplane  bundle of $\C P^n$.
An orthonormal basis of the space $H^0(M,L^m)$ can be
represented by
\begin{equation}\label{6-1}
\sqrt{\frac{(m+n)!}{P!}}z^P
\end{equation}
for multiple index $P\in\Z^n_+$ with $|P|=m$, where
$P!=p_1!\cdots p_n!$ for $P=(p_1,\cdots,p_n)$.

 We
shall first compute concretely the finite 
dimensional vector space $V$ in
Theorem~\ref{fud}.

Consider the open set $U_0$ of $\cp$  where the local
coordinate is $(z_1,\cdots,z_n)$ and the homogeneous
coordinate is represented by $[1,z_1,\cdots,z_n]$.
Since ${\mathbb C}P^n$ is a symmetric space, we only 
need to consider
the expansion at the point $x_0=[1,0,\cdots,0]$. The local 
coordinate of $x_0$ is $(0,\cdots,0)$. Thus in the following
we  sometimes use $0$ to represent the point  $x_0$.

 Let
$S_0=\sqrt{\frac{(m+n)!}{m!}}$ be the section under the
standard local trivialization of $H$ on $U_0$. The Hermitian
metric on
$L$ is defined by
$h=1/(1+|z|^2)$, and the K\"ahler metric is defined by
$\omega=\bb\pa\bar\pa\log (1+|z|^2)$. The pointwise norm
of the section $S_0$ at $x_0=(0,\cdots,0)$ is
\begin{equation}\label{6-2}
||S_0||^2={\frac{(m+n)!}{m!}}\cdot\frac{1}{(1+|z|^2)^m}.
\end{equation}
 Under the
basis 
~\eqref{6-1}, 
$I_{0\alpha}=0$, $\sigma(0)=\frac{(m+n)!}{m!}$ and 
$I_{00}^{-1}=1$. Using~\eqref{prop51}, we have
\[
\left.\frac{d}{dt}\right|_{t=0}\sigma(t)
=-\frac{(m+n)!}{m!n!}\int_M(m\phi-\Delta\phi)||S_0||^2\omega^n.
\]
Substituting ~\eqref{6-2} into thte above equation, we have
\begin{equation}\label{6-9}
\left.\frac{d}{dt}\right|_{t=0}\sigma(t)=
-\frac{1}{\pi^nn!}\left(\frac{(m+n)!}{m!}\right)^2\int_{\C^n}(m\phi-\Delta\phi)
\frac{1}{(1+|z|^2)^{m+n+1}}dV_0,
\end{equation}
where $dV_0$ is the Euclidean volume form  of $\C^n$.

The following identity is elementary and will be used
repeatedly:
\begin{equation}\label{ele1}
\int_{\C^n}\frac{|z^P|^2}{(1+|z|^2)^{m+n+1}} dV_0
=\pi^n\frac{P!(m-|P|)!}{(m+n)!},
\end{equation}
where $|P|\leq m$.

\smallskip

\begin{lemma}\label{lem61}
There is an asymptotic expansion of the right hand
side of ~\eqref{6-9}:
\[
\left.\frac{d}{dt}\right|_{t=0}\sigma(t)\sim
m^n(\xi_0+\frac{\xi_1}{m}+\cdots)
\]
at the point $x_0$
where $\xi_i$ can be represented as 
$\xi_i=f_i(\Delta_c)\phi(0)$, $i\geq 1$ for
 polynomials $f_i$. The operator
$\Delta_c$ is defined in ~\eqref{3-1-1}.
\end{lemma}

{\bf Proof.} 
The existence of the expansion is from Theorem~\ref{thm52}.
Assuming $\ph(0)=0$, the Taylor expansion of
$\phi$ at $x_0$ is
\[
\ph\sim\sum_{|P|+|Q|>0}\frac{1}{P!Q!}a(P,Q)z^P\bar z^Q.
\]
Using~\eqref{ele1}, we have
\begin{equation}\label{6-8}
\int_{\C^n}\frac{1}{(1+|z|^2)^{m+n+1}}\ph dV_0
\sim\pi^n\sum_P\frac{(m-|P|)!}{P!(m+n)!}
a(P,P).
\end{equation}
We also have
\[
\Delta_c^k\ph(0)=\sum_{|P|=k}\frac{k!}{P!}a(P,P).
\]
From the above equation, ~\eqref{6-8} becomes
\[
\int_{\C^n}\frac{1}{(1+|z|^2)^{m+n+1}}\ph dV_0
\sim\pi^n\sum_{k=1}^\infty\frac{(m-k)!}{k!(m+n)!}
\Delta_c^{k}\ph(0).
\]
We thus have the expansion
\begin{equation}\label{6-9-1}
\int_{\C^n}\frac{1}{(1+|z|^2)^{m+n+1}}\ph dV_0 
\sim\frac{1}{m^{n+1}}(\eta_0+\frac{\eta_1}{m}+\cdots),
\end{equation}
where the coefficients are all polynomials of
$\Delta_c$ acting on $\phi$ at $0$. The lemma follows from
~\eqref{6-9-1}.

\qed

The following proposition is purely combinatoric:

\begin{prop}\label{comb}
There are polynomials
\[
f_k(t)=\sum_{l=0}^k a_{k,l}t^l
\]
of degree $k$ such that 
\begin{equation}\label{6-1-1}
\left\{
\begin{array}{l}
a_{k,0}=0,\\
a_{k,k}=1,\\
a_{k,k+1}=0,\\
\Delta^k\phi(0)=f_k(\Delta_c)\phi (0),
\end{array}
\right.
\end{equation}
where $\phi$ is a smooth function,
$\Delta$ is the Laplacian of $\cp$, and $k\in\N$.
\end{prop} 

{\bf Proof.} If $k=1$, then we choose $f_k(t)=t$. 
~\eqref{6-1-1} is valid.
Using the
mathematical induction,  we
assume that for $k\geq 1$,
\[
\Delta^k\phi(0)=\sum_{l=0}^ka_{k,l}\Delta_c^l\phi(0)
\]
for constants $a_{k,l}\,(0\leq l\leq k)$.
We wish to construct constants $a_{k+1,l}$ with $0\leq l\leq
k+1$ such that ~\eqref{6-1-1} is true for $k+1$.

We need the following lemma:
\begin{lemma}\label{lem62}
Define
\begin{equation}\label{6-4}
\left\{
\begin{array}{l}
a_{k+1,0}=0,\\
a_{k+1,k+1}=1,\\
a_{k+1,k+2}=0,\\
a_{k+1,l}=a_{k,l-1}+l(2l+n-1)a_{k,l}
+l^2(l+1)(l+n)a_{k,l+1},\quad 0<l<k+1.
\end{array}
\right.
\end{equation}
Then we have
\begin{equation}\label{6-5}
\Delta^{k+1} |z^P|^2(0)=\sum_{l=0}^{k+1}
a_{k+1,l}\,\Delta_c^l|z^P|^2(0)
\end{equation}
for $|P|\leq k+1$.
\end{lemma}

{\bf Proof.} Firstly, if $|P|=k+1$, then 
\[
\Delta^{k+1}|z^P|^2(0)=\Delta_c^{k+1}|z^P|^2(0)
\]
and
\[
\Delta_c^l|z^P|^2(0)=0
\]
for $l<k+1$. 
Thus in this case, ~\eqref{6-5} holds true.
Note that 
\[
\Delta=(1+|z|^2)(\delta_{ij}+z_i\bar z_j)\frac{\pa^2}{\pa
z_i\pa \bar z_j}.
\]
Now let $|P|=l<k+1$. Then we have

\begin{equation}\label{6-10}
\Delta|z^P|^2=\sum_ip_i^2|z^{P_i}|^2+l^2|z^P|^2+\sum_{i,j}
p_i^2|z^{Q_{ij}}|^2+\sum_{i}l^2|z^{R_i}|^2,
\end{equation}
where $P_i, Q_{ij}$ and $R_i$ are defined as follows:
let $e_j=(0,\cdots,\underset{j}{1},\cdots,0)$ for $1\leq
j\leq n$. Recall that a multiple index $S\geq 0$ iff
all of its components $s_i\geq 0$. For $1\leq i,j\leq n$,
define
\begin{enumerate}
\item $P_i=P-e_i$ if $P_i\geq 0$, otherwise $P_i=P$;
\item $Q_{ij}=P_i+e_j$;
\item $R_i=P+e_i$.
\end{enumerate}
The lemma follows from ~\eqref{6-10},
the math induction, and the following
identity
\begin{equation}\label{6-11}
\Delta_c ^l|z^P|^2(0)=l!P!.
\end{equation}

\qed

\begin{lemma}\label{lem63}
If $P\neq Q$, then 
\begin{equation}\label{6-12}
\Delta^k z^P\bar z^Q(0)=0.
\end{equation}
\end{lemma}

{\bf Proof.} We use the mathematical induction again. If
$k=1$, the theorem is true. Thus we assume that ~\eqref{6-12}
is true for 
any $k\leq s$. A straightforward computation gives
\[
\Delta z^P\bar z^Q=\sum
a_{R,S}z^R\bar z^S,
\]
where $R\neq S$. Thus we have
\[
\Delta^{s+1}z^P\bar z^Q(0)=
\sum a_{R,S}\Delta^sz^R\bar z^S(0)=0.
\]
The lemma is proved.

\qed

{\bf Continuation of the proof of Proposition~\ref{comb}.} By
linearity, we just need to verify
\begin{equation}\label{6-13}
\Delta^k\ph(0)=f_k(\Delta_c)\ph(0)
\end{equation}
when $\ph=z^P\bar z^Q$. If $|P|+|Q|>2k$, then both sides
of the above equation are zero. If $P\neq Q$, then by
Lemma~\ref{lem63}, both sides of the above equation are
zero. If $P=Q$, then by Lemma~\ref{lem62}, ~\eqref{6-13} holds
true.

\qed

Using Lemma~\ref{lem61}, Proposition~\ref{comb} and the 
homogenity of $\cp$, we
have the following

\begin{theorem}
Let $x\in\cp$. Then
in the expansion
\[
\left.\frac{d}{dt}\sigma(t)\right|_{t=0}\sim
m^n(b_0+\frac{b_1}{m}+\cdots),
\]
the coefficients can be represented by
\[
b_i(x)=f_i(\Delta)\ph(x),\quad i\geq 0,
\]
where $\Delta$ is the Laplacian of  $\cp$.
\end{theorem}

{\bf Proof.} Let $x=x_0=[1,0,\cdots,0]$. Then by Lemma
~\ref{lem61}, $b_i$, $i\geq 1$ can be represented as 
$g_i(\Delta_c)$ for polynomials $g_i$. Using
Proposition~\ref{comb}, we can write
$b_i=f_i(\Delta)\phi(x_0)$. The general case follows from 
the homogeneity of $\cp$.

\qed

If the log terms of the Szeg\"o kernel  with respect to 
the metrics $h_t=he^{t\phi}$ are zero for small $t$, then 
$b_i(x)\equiv 0$ for $i>n$ by
Theorem~\ref{thm23} and Theorem~\ref{34}. In particular
$b_{n+1}=f_{n+1}(\Delta)\phi=0$.

\begin{lemma}
If $f_{n+1}(\Delta)\ph=0$, then $\ph$ is zero or is an
eigenfunction of $\Delta$.
\end{lemma}

{\bf Proof.} We assume that
\[
f_{n+1}(t)=\prod_{i=0}^s(t-\mu_i),
\]
where $\mu_i$ are complex numbers. Then we have
\[
\prod_{i=0}^s(\Delta-\mu_iI)\ph=0.
\]

Define
\[
\psi_i=\prod_{k=i}^s(\Delta-\mu_kI)\ph
\]
for $0\leq i\leq s$. Then we have $\psi_0=0$,
and
\[
(\Delta-\mu_{i-1}I)\psi_i=\psi_{i-1},\quad i\geq 1.
\]
Assuming that
$\Delta\psi_{i-1}=-\lambda\psi_{i-1}$ for some $i$, then if
$\lambda+\mu_i\neq 0$,  we have
\[
\psi_i=-\frac{1}{\lambda+\mu_{i-1}}\psi_{i-1}.
\]
In particular, $\Delta\psi_i=-\lambda\psi_i$. If
$\lambda+\mu_{i-1}=0$, then we still have
$\Delta\phi_i=-\lambda\phi_i$. At this time,
$\phi_{i-1}=(\Delta-\mu_{i-1})\psi_i=0$.

In this case, $\psi_{i-1}=0$ and
$\Delta\psi_i=-\lambda\psi_i$. Using the mathematical
induction,
we have $\Delta\ph_{s+1}=-\lambda\ph_{s+1}$ for some
nonnegative real number $\lambda$, where $\ph_{s+1}=\ph$.

\qed

We wish to prove that $\lambda=-(n+1)$. To this purpose,
we define
\begin{equation}\label{def6}
\ph_k=\frac{1}{(1+|z|^2)^k}
\end{equation}
for $k=0,1,2,\cdots$. A straightforward computation
gives
\begin{equation}\label{6-17}
\Delta\ph_k=-k(k+n)\ph_k+k^2\ph_{k-1}
\end{equation}
for $k\geq 1$. Using the integration by parts, we have
\begin{equation}\label{6-18}
(k(k+n)-\lambda)\int_{\cp}\ph\ph_k=k^2\int_{\cp}\ph\ph_{k-1}.
\end{equation}
If $\lambda\neq k(k+n)$ for any integer $k$, then since
\[
\int_{\cp}\ph=0,
\]
using ~\eqref{6-18}, we see that
\[
\int_{\cp}\ph\ph_k=0
\]
for any $k$. By ~\eqref{6-9}, this implies that
\[
\left.\frac{d}{dt}\sigma(t)\right|_{t=0}=0.
\]
However, this is not possible, because we have
\[
\left.\frac{d}{dt}\right|_{t=0}a_1=\Delta(\Delta+(n+1)I))\ph=0,
\]
which implies that $\lambda=n+1$. Now assume that
$\lambda=k_0(k_0+n)$~\footnote{The proof also implies that
the eigenvalues of $\cp$ must be of the form $k(k+n)$ for
 $k\in\Z$.}. Then by~\eqref{6-18}, we have
\[
\int_{\cp}
\ph\ph_m=\frac{(m!)^2/(k_0!)^2}{(m-k_0)!(m+k_0+n)!/(2k_0+n)!}
\int_{\cp}\ph\ph_{k_0}.
\]
Thus ~\eqref{6-9} becomes
\[
\left.\frac{d}{dt}\right|_{t=0}\sigma(t)
={\rm const}\cdot 
\frac{(m+n)\cdots(m-k_0+1)(m+k_0(k_0+n))}
{(m+k_0+n)\cdots (m+n+1)}.
\]
The above expression is a polynomial of $m$ if and only if
$k_0=1$. Thus we have proved

\begin{theorem}\label{lamb}
For the complex projective space $\cp$, the vector space $V$
in Theorem~\ref{fud} is contained in the eigenspace of the
first eigenvalue of
$\cp$.
\end{theorem}

\qed

\section{The proof of the main theorem} 
%2003.2.18.9.53
In the last section, we proved that the 
vector space $V$ in Theorem~\ref{fud}
 is contained in the eigenspace of the eigenvalue $(n+1)$. 
Among all the eigenspaces, the eigenspace of the
eigenvalue $(n+1)$ is special. We have the following 
well known result (cf. ~\cite{CT1}):

\begin{lemma}\label{wellknown}
We use the same notations as in the previous section. Then
the first eigenvalue of $\cp$ with the Fubini-Study metric 
is $(n+1)$.
Let $\ph$ be an eigenfunction of the eigenvalue $(n+1)$.
Let the $(1,0)$ vector field  $X$ on $\cp$ be defined by
\[
X=g^{i\bar j}\ph_{\bar j}\frac{\pa}{\pa z^i}.
\]
Then $X$ is holomorphic.
\end{lemma}

\qed

The automorphism of $\cp$ can be represented by a
nonsingular
$(n+1)\times (n+1)$ matrix $a_{ij}$. That is, for
any such matrix, the linear map
\begin{equation}\label{easy}
f: \cp\rightarrow \cp,\quad Z_i\mapsto\sum_j a_{ij}Z_j
\end{equation}
defines an automorphism. The Bergman potential is invariant
under the automorphism, so are the coefficients 
in the Tian-Yau-Zelditch expansion.  Thus in order to prove
the theorem, we must get rid of these automorphisms.
The method below is similar  to that in Bando and Mabuchi
~\cite{BM} in the study of the uniqueness of 
K\"ahler-Einstein metrics on Fano manifolds. However,
we use the notations and the results in~\cite{CT2}.

Let $\omega_0$ be the \ka form of the standard Fubini-Study
metric. Let $\omega_\rho$ be defined as
\[
\omega_\rho=\omega_0+\bb\pa\bar\pa\rho=
f_\rho^*\omega_0,
\]
where $\rho$ 
is  a real valued function defined by
\[
\rho=\log\frac{\sum_i|\sum_j a_{ij}Z_j|^2}{\sum_i|Z_i|^2}
\]
for a nonsingular
 $(n+1)\times (n+1)$ matrix $(a_{ij})$. $f_\rho$ is the
automorphism defined by $\rho$ as in ~\eqref{easy}.

\begin{definition}
Any \ka form $\omega_\ph=\omega_0+\bb\pa\bar\pa\ph$ is
called centrally positioned with respect to the 
metric $\omega_\rho=\omega_0+\bb
\pa\bar\pa\rho$ (which is the
Fubini-Study metric), if
\[
\int_{\cp}(\rho+f^*_\rho(\ph))f_\rho^*(\theta)\omega_\rho^n
=0
\]
for any $\theta\in \Lambda_{n+1}(\omega_0)$, where
$\Lambda_{n+1}(\omega_0)$ is the space of eigenfunctions
with respect to
the first eigenvalue $(n+1)$ of $\cp$.
\end{definition}

For any $\rho$, the \ka metric defined by  the function
$\rho+f_\rho^*(\ph)$ and $\ph$ differ only by 
the automorphism $f_\rho$. So they are essentially 
equivalent.
In particular, the vanishing of the log term of the Szeg\"o
kernel for one metric implies the vanishing of the log
term of the Szeg\"o
kernel for the other. We thus need to choose
the best representative among all these equivalent metrics.

The following proposition 
shows that the best representative always exists.

\begin{prop}~\cite{CT2}
Using the notations as above, then for any $\ph$, there is 
a function $\rho$ such that $\omega_\ph$  is centrally
positioned.
\end{prop}

\qed

For our purpose, we only  need the existence of $\rho$
when $\ph$ is small. However, we need to estimate the 
solution $\rho$. We will use the following 
method of contraction principal to construct the 
solution. 

\begin{lemma}\label{lem72}
Using the notations as above, then for any $\eps>0$, there
is an $\eta$ such that if $||\ph||_{C^0}<\eta$
then $||\rho||\leq\eps$.
\end{lemma}

By  Lemma~\ref{wellknown}, we  know that the dimension of the
eigenspace of the first eigenvalue is equal to the dimension
of the space of  holomorphic vector fields, which is
equal to $(n+1)^2-1$, the dimension of the automorphism 
group of $\cp$.

Let ${\frk p}$ be the space of all $(n+1)\times (n+1)$
 traceless  Hermitian matrices. The real dimension of
${\frk p}$ is $(n+1)^2-1$. Let $\theta_1,\cdots,\theta_s$
$(s=(n+1)^2-1)$ be a (real orthonormal basis of 
$\Lambda_{n+1}(\omega_0)$, the eigenspace of the eigenvalue
$(n+1)$
of the Fubini-Study metric. By direct calculation, we see
that
\[
\left\{
\begin{array}{ll}
\frac{Z_i\bar Z_j}{\sum_k|Z_k|^2}, & 0\leq i, j\leq n,\\
\frac{|Z_i|^2-|Z_0|^2}{\sum_k|Z_k|^2}, & 1\leq i\leq n,
\end{array}
\right.
\]
are a (complex) basis of eigenfunctions. For any
$(n+1)\times (n+1)$ matrix $A\in {\frk p}$. Consider
the real eigenfunction
\[
\frac{\sum_{i,j}a_{ij}Z_i\bar Z_j}{\sum_k|Z_k|^2}.
\]
If we represented the function as
\[
 \frac{\sum_{i,j}a_{ij}Z_i\bar Z_j}{\sum_k|Z_k|^2}
=\sum_{i=1}^s b_i\theta_i,
\]
then we can define
\[
L: {\frk p}\rightarrow \R^s,\quad
L(A)=(b_1,\cdots,b_s).
\]
Since a real basis of $\Lambda_{n+1}(\omega_0)$
can be represented by $(n+1)^2-1$ functions
\[
{\rm Re}\,\frac{Z_i\bar Z_j}{\sum_k|Z_k|^2},\,\,
{\rm Im}\,\frac{Z_i\bar Z_j}{\sum_k|Z_k|^2},\,\,
\frac{|Z_i|^2-|Z_0|^2}{\sum_k|Z_k|^2},
\]
$L$ is an invertible linear map.

Let $U$ be a small neighborhood of ${\frk p}$ at origin.
Let $A\in{\frk p}$. Let $e^A=(\tilde a_{ij})$. Define
\[
\rho_A=\log\frac{\sum_i|\sum_j\tilde a_{ij}Z_j|^2}
{\sum_k|Z_k|^2}.
\]
Let $\theta=(\theta_1,\cdots,\theta_s)$. Define the following 
nonlinear operator
\begin{equation}\label{def5}
T: U\times C^0(M)\rightarrow {\frk p}:
(A,\ph)\mapsto A-\frac 12L^{-1}\int_{\cp}
(\rho_A+f^*_{\rho_A}(\ph))f^*_{\rho_A}(\theta)
\omega^n_{\rho_A}.
\end{equation}

We shall see that there is a fixed point of the operator
if the $||\ph||_{C^0}$ is fixed and is small. Let
$T(A)=T(A,\ph)$. Then there
is a constant
$C$, such that
\[
||T(0)||\leq C||\ph||_{C^0}.  
\]
We also have
\[
||T(B)-T(A)||
=
||B-A-\frac{L^{-1}}{2}\int_{\cp}((f^{-1}_{\rho_A})^*(\rho_A))
-(f^{-1}_{\rho_B})^*(\rho_B))\theta\omega_0^n||.
\]
If $A, B$ are small, then there is a constant $K$, such that
\[
||T(B)-T(A)||\leq K||B-A||^2.
\]
If $||A||, ||B||\leq 1/4K$, then we have
\[
||T(B)-T(A)||\leq\frac 12 ||B-A||.
\]
We choose $||\ph||_{C^0}$ to be small enough so that
\[
||T(0)||\leq{\rm Min}(\frac\eps 2, \frac{1}{8K}).
\]
Then
\[
||T^k(0)-T^{k-1}(0)||\leq\frac{1}{2^{k-1}}{\rm Min}
(\frac \eps 2, \frac{1}{8K})
\]
for $k=1,2,\cdots$.
Thus we have 
$
\sum_{k=1}^\infty ||T^k(0)-T^{k-1}(0)||\leq\eps,
$
from which we have
\[
\underset{k\rightarrow\infty}{\lim} T^k(0)=A
\]
exists, $||A||\leq \eps$ and $A$ is a fixed point of $T$.
The lemma follows from the fact that $\omega_\phi$ is
centrally positioned with respect to $\rho_A$.

\qed

{\bf Proof of the Theorem~\ref{fud1}.} We assume that there is
a sequence $\ph_i$, $i\geq 1$ such that
$||\phi_i||_{C^{2n+4}}\rightarrow 0$, and 
\begin{equation}\label{contra}
a_{n+1}(x,he^{\ph_i})=0
\end{equation}
for $i\geq 1$.
By the above lemma, we can replace $\ph_i$ by
$\tilde\ph_i=\rho_i-f_{\rho_i}^*(\ph_i)$ and we still
have
\[
a_{n+1}(x,he^{-\rho_i+f_{\rho_i}^*(\ph_i)})=0
\]
with
\begin{equation}\label{7-7}
\int_{\cp}(\rho_i-f_{\rho_i}^*(\ph_i))f^*_{\rho_i}\theta
\omega_{\rho_i}^n=0.
\end{equation}
By Lemma~\ref{lem72}, we have
\begin{equation}\label{7-8}
||\rho_i-f_{\rho_i}^*(\ph_i)||_{C^{2n+4}}=\eps_i\rightarrow
0.
\end{equation}
By Corollary~\ref{cor31}, we see that there is a constant $C$
such that
\[
||\rho_i-f_{\rho_i}^*(\ph_i)||_{C^{2n+6}}\leq C\eps_i
\]
Thus there is a subsequence of $\xi_i=
(f_{\rho_i}^*(\ph_i)-\rho_i)/\eps_i$,
for which we still denote
as
$\xi_i$, converges to some $\xi\neq 0$ in the $C^{2n+4}$
norm.  Furthermore, using Lemma~\ref{lem72} again, we get
\begin{equation}\label{7-11}
\int_{\cp}\xi\theta\omega_0^n=0.
\end{equation}
On the other side, by ~\eqref{f-1}, we have
\begin{equation}\label{7-12-1}
\left|a_{n+1}(x,he^{\eps_i\xi})-a_{n+1}(x,h)-\left.
\frac{d}{ds}\right|_{s=0}a_{n+1}(x,he^{s\xi})\eps_i\right|
\leq\eps_i^2.
\end{equation}
Using~\eqref{f-2},  we have
\begin{equation}\label{7-12-2}
\left|a_{n+1}(x,he^{-\rho_i+f^*_{\rho_i}(\ph_i)})-a_{n+1}(x,he^{\eps_i\xi})
\right|\leq
C\eps_i||\frac{\rho_i-f^*_{\rho_i}(\ph_i)}
{\eps_i}-\xi||_{C^{2n+4}}.
\end{equation}
By assumption, $a_{n+1}(x,he^{-\rho_i+f^*_{\rho_i}(\ph_i)})
=0$. Since $h$ is the standard metric of $\cp$,
$a_{n+1}(x,h)=0$. Thus from~\eqref{7-12-1} and~\eqref{7-12-2},
we have
\[
\left.\frac{d}{ds}\right|_{s=0}a_{n+1}(x,he^{s\xi})=0.
\]
By Theorem~\ref{lamb}, this implies
\[
\Delta\xi=-(n+1)\xi.
\]
From~\eqref{7-11}, $\xi\equiv 0$. This is a contradiction.
So for $||\phi_i||_{C^{2n+4}}$ small, $a_{n+1}(x,
he^{\phi_i})\neq 0$. This proves the theorem.

\qed

\bibliographystyle{abbrv}
%\bibliography{::::bib:bib}
\bibliography{bib}

\end{document}